\def\obs#1{{\bf (*** #1 ***)} }

 \def\obs#1{}     
\documentclass[twoside,letterpaper,draft,11pt]{amsart}
\usepackage{amsmath,amsthm,latexsym,xspace,amscd,amssymb,mathtools,color,enumerate}

\usepackage[mathscr]{eucal}
\usepackage{amssymb}
\usepackage[all]{xy}
\makeatletter\renewcommand\theenumi{\@roman\c@enumi}\makeatother

\newtheorem{teo1}{Theorem}[section]

\newtheorem{lem1}[teo1]{Lemma}
\newtheorem{cor1}[teo1]{Corollary}
\newtheorem{prop1}[teo1]{Proposition}

\newtheorem{exe}[teo1]{Example}
\newtheorem{remark}[teo1]{Remark}

\newcommand{\m}{{}^{-1}}

\newcommand{\vu}{\vspace{.1cm}}
\newcommand{\vd}{\vspace{.2cm}}

\newcommand{\Res}{\mathtt{Res}}
\newcommand{\Ext}{\mathtt{Ext}}
\newcommand{\Par}{\textsf{Par}_\G(\A)}
\newcommand{\GT}{\textsf{GrT}_\G(\A)}
\newcommand{\Dat}{\textsf{D}_{\G,\tau(x)}(\A)}
\newcommand{\NDat}{\textsf{D}_{\G,\lambda(z)}(\A)}
\newcommand{\DA}{\textsf{D}_\G}
\newcommand{\GDA}{\textsf{C}_\G}
\newcommand{\nres}{\mathtt{R}}
\newcommand{\next}{\mathtt{E}}

\newcommand{\id}{{\rm id}}

\newcommand{\I}{{I}}

\newcommand{\A}{{A}}

\newcommand{\af}{\alpha}
\newcommand{\bt}{\beta}
\newcommand{\lb}{\lambda}

\newcommand{\gm}{\gamma}

\newcommand{\ta}{\tau}

\newcommand{\somai}{\sum_{1\leq i\leq r}}
\newcommand{\somay}{\sum_{y\in \G_0}}
\newcommand{\somag}{\sum_{g\in \G}}
\newcommand{\somaz}{\sum_{z\in\G_0}}
\newcommand{\somal}{\sum_{l\in\G(x)}}

\newcommand{\G}{\mathcal{G}}

\def\ndv{\ {\mid \kern -0.7 em {\scriptstyle \not}} \ \ }

\def\nd{\ {\mid \kern -0.4 em {\scriptstyle \not}} \ \ }

\begin{document}

\thispagestyle{empty}

\title[Restriction and extension of partial actions]{Restriction and extension of partial actions}

\author[D. Bagio]{Dirceu Bagio}
\address{ Departamento de Matem\'atica, Universidade Federal de Santa Maria, 97105-900\\
	Santa Maria-RS, Brasil}
\email{ bagio@smail.ufsm.br}

\author[A. Paques]{Antonio Paques }
\address{Instituto de Matem\'atica e Estat\'istica, Universidade federal de Porto Alegre, 91509-900\\
	Porto Alegre-RS, Brazil}
\email{paques@mat.ufrgs.br}

\author[H. Pinedo]{H\'ector Pinedo}
\address{Escuela de Matematicas, Universidad Industrial de Santander, Cra. 27 Calle 9  UIS
	Edificio 45, Bucaramanga, Colombia}
\email{ hpinedot@uis.edu.co}


\date{\today}
\begin{abstract} Given a partial action $\af=(\A_g,\af_g)_{g\in \G}$ of a connected groupoid $\G$ on a ring $A$ and an object $x$ of $\G$, the isotropy group $\G(x)$ acts partially on the ideal $A_x$ of $A$ by the restriction of $\alpha$. In this paper we investigate the following reverse question: under which conditions a partial group action of $\G(x)$ on an ideal of $A$ can be extended to a partial groupoid action of $\G$ on $A$? The globalization problem and some applications to the Morita and Galois theories are also considered, as extensions of similar results from the group actions case.
\end{abstract}

\maketitle
\setcounter{tocdepth}{1}

\tableofcontents

\section{Introduction}
The notion of partial group action, which was introduced  in
\cite{EX},  is  an efficient tool to study the structure of  C*-algebras (see \cite{EX2}).   Partial group actions on sets can 
be described in terms of premorphisms. 
In fact, in \cite{KL} it was shown that a  partial action of a group $G$ on a set $X$ is a unital premorphism from $G$ to the inverse monoid  $\mathcal{I}(X)$ of partial bijections of $X$. In the same way, N. D. Gilbert introduced in \cite{Gi} the concept of partial actions of ordered groupoids on sets.  
Motivated by N. D. Gilbert's work, the authors of \cite{BFP} defined partial actions of ordered groupoids on rings. After that, partial actions of groupoids on rings was considered by
D. Bagio and A. Paques in \cite{BP}. From then on, several questions involving partial groupoid actions were investigated. For instance, ring-theoretical properties of the partial skew groupoid ring were studied in \cite{BPi} and \cite{NYOP} and the Galois correspondence theorem for groupoid actions was given in \cite{PT}.

It is well-known that every groupoid is a disjoint union of its connected components and that a partial action of a groupoid on a ring is uniquely determined by the respective partial actions of its connected components (see \cite[Proposition 2.3]{BPi}). Therefore, we can reduce the study of partial groupoid actions to the connected case. \vspace{0.07cm}

Let $\af=(\A_g,\af_g)_{g\in \G}$ be a partial action of a connected groupoid $\G$ on a ring $A$ and $x$ an object of $\G$. We highlight two important facts. First, the isotropy group $\G(x)$ acts partially on the ideal $A_x$ of $A$ by the restriction of $\alpha$. Second, the groupoid $\G$ is isomorphic to $\G_0^2\times \G(x)$, where $\G_0$ is the set of objects of $\G$ and $\G_0^2$ is the coarse groupoid associated to $\G_0$. Therefore it is natural to expect that, under appropriate conditions, we can extend a partial action of $\G(x)$ on an ideal of $A$ to a partial action of $\G$ on $A$. The main purpose of this paper is to determine such conditions as well as to construct such extensions. \vspace{0.07cm}

We organize the paper as follows. The background about groupoids and partial actions that will be used later is given in Section 2. The main results are in Section 3. We present sufficient conditions to extend partial group actions to partial groupoid actions in $\S\,3.1$. The explicit construction of the extensions is given by \eqref{afcon} and \eqref{dom} in $\S\,3.2$. Some examples that illustrate the construction are in $\S\,3.3$. The globalization problem of partial groupoid actions obtained in this setting is discussed in Section 4.
Applications to the Morita and Galois theories are given in Section 5.

\subsection*{Conventions}\label{subsec:conv}
Throughout this work, the rings are associative and not necessarily unital. Given a groupoid $\G$ and a morphism $g:x\to y$ in $\G$, we will often write $g\in \G$. We will work with a functor that extends partial group actions. This functor will be denoted by $\Ext$ although it is not related to the extension functor used in the homological algebra.

\section{Preliminaries}
In this section we recall the background about groupoids and partial actions that will be used in the work.

\subsection{Groupoids}
A {\it groupoid} $\G$ is a small category in which every morphism is an isomorphism. The set of the objects of $\G$ will be denoted by $\G_0$. As usual, each object $x$ of $\G$ is identified with its identity morphism $\id_x$, whence $\G_0$ can be seen also as the set of the  identity elements of $\G$. Given $g\in \G$, the {\it source} (or domain) and the {\it  target } (or codomain) of $g$ will be denoted by $s(g)$ and $t(g)$, respectively.

\vu

The composition of morphisms of a groupoid $\G$ will be denoted via concatenation and we will compose from right to left. We will write $\G_s\times\!_t\G$ to denote the set of ordered pairs of $\G$ which are composable, that is, $\G_s\times\!_t\G=\{(g,h)\in\G\times\G \ :\ s(g)=t(h)\}$. Notice that
\[ s(x)=x=t(x),\quad  s(g)=g^{-1}g,\quad t(g)=gg^{-1},\quad s(gh)=s(h)\quad \text{and} \quad t(gh)=t(g),\]
for all $x\in \G_0$, $g\in \G$ and $(g,h)\in\G_s\!\times\!\!\,_t\G$. These relations will be freely used along the text. For $x,y\in\G_0$ we set $\G(x,y)=\{g\in \mathcal{G}: s(g)=x\,\,\text{and}\,\,t(g)=y\}$. Particularly $\G(x):=\G(x,x)$ is a group, called the \emph{isotropy group associated to $x$.}

A groupoid $\G$ is called \emph{connected (or Brandt groupoid)} if $\mathcal{\G}(x,y)\neq \emptyset$ for all $x,y\in \G_0$. It is well-known that any groupoid is a disjoint union of connected subgroupoids. In fact, we have the following equivalence relation on $\G_0$: for all $x,y\in \G_0$,
\[x\sim y \,\,\text{ if and only if }\,\, \G(x,y)\neq \emptyset.\]
Every equivalence class $X\in\G_0/\!\!\sim$ determines a subgroupoid $\G_X$ of $\G$ in the following way. The set of the objects of $\G_X$ is $X$ and $\G_X(x,y)=\G(x,y)$ for all $x,y\in X$. By construction, $\G_X$ is a full connected subgroupoid of $\G$ which is called the {\it connected component of $\G$ associated to $X$}. Clearly, $\G=\dot\cup_{X\in \G_0/\!\sim}\G_X$, that is, $\G$ is the disjoint union of its connected components.

The structure of a connected groupoid is well-known. It is the cartesian product of a coarse groupoid by a group (see Proposition \ref{group:connec} below). We recall that the {\it coarse groupoid} associated to a nonempty set $Y$ is the groupoid $Y^2=Y\times Y$, whose set of objects is $Y$ and  $\{(y,y)\ |\ y\in Y\}$ is the set of its identities. The composition in $Y^2$ is given by: $(z,w)(y,z)=(y,w)$, and it is immediate to see that $(y,z)^{-1}=(z,y)$, $s(y,z)=y$ and $t(y,z)=z$, for all $y,z,w\in Y$.

\begin{prop1}\label{group:connec}
	Let $\G$ be a connected groupoid, $x\in \G_0$, $\tau_y\in \G(x,y)$ for all $y\in \G_0$, $y\neq x$, and $\tau_x=x$. Then $\psi:\G\to \G_0^2\times \G(x)$ given by
	\[\psi(g)=((s(g),t(g)), g_x),\quad \text{where}\,\,\, \text g_x:=\ta^{-1}_{t(g)}g\ta_{s(g)},\] is an isomorphism of groupoids.
\end{prop1}
\begin{proof}
Let $(g,h)\in \G_s\times\!_t\G$. Since $s(gh)=s(h)$, $t(gh)=t(g)$ and $(gh)_x=g_xh_x$ it follows that $\varphi$ is a groupoid morphism. Let $g\in \G$ and suppose that $\varphi(g)=((y,y),x)$ for some $y\in\G_0$, i.~e. $\varphi(g)$ is an object of $\G_0^2\times \G(x)$. Thus $s(g)=t(g)=y$ and $x=g_x=\ta_y^{-1}g\ta_y$ which implies that $g=\ta_y\ta_y^{-1}=x$. This ensures that $\varphi$ is injective. Given an element $((y,z),h)\in \G_0^2\times \G(x)$, consider $g=\tau_zh\tau\m_y\in \G$. Notice that $g_x=h$ and whence $\varphi(g)=((y,z),h)$, that is,  $\varphi$ is surjective, so an isomorphism of groupoids. \end{proof}

\subsection{Partial and global actions of groupoids}
We start by recalling from \cite{BP} the definition of partial groupoid action on rings. A \emph{partial action} of a groupoid $\G$ on  a ring $\A$ is a pair
$\af=(\A_g,\af_g)_{g\in \G}$ that satisfies
\begin{enumerate}[\hspace{.35cm}(P1)]
	\item for each $g\in \G$,  $\A_{t(g)}$ is an ideal of $\A$, $\A_g$ is an ideal of $\A_{t(g)}$ and $\af_g:\A_{g\m}\to \A_g$ is an isomorphism of rings,\smallbreak
	\item for each $x\in \G_0$, $\af_x=\id_{A_x}$ is the identity map  of $\A_x$, \smallbreak
	\item $\af_h^{-1}(\A_{g^{-1}}\cap\A_h)\subseteq A_{(gh)^{-1}}$, for all $(g,h)\in \G_s\!\times\!\!\,_t\G $,\smallbreak
	\item $\af_g(\af_h(a))=\af_{gh}(a)$, for all $a\in \af_h^{-1}(\A_{g^{-1}}\cap\A_h)$ and $(g,h)\in \G_s\!\times\!\!\,_t\G $.
\end{enumerate}

\begin{remark}\label{obs-inicial}
{\rm It follows from  (P3) and (P4) that $\alpha_{g}\alpha_{h}\leq \alpha_{gh}$, that is, $\af_{gh}$ is an extension of $\af_g\af_h$, for all  $(g,h)\in \G_s\!\times\!\!\,_t\G $.} 	
\end{remark}

The partial action $\af$ is said to be {\it global} if $\af_g\af_h=\af_{gh}$, for all $(g,h)\in \G_s\!\times\!\!\,_t\G$. By Lemma 1.1 of \cite{BP}, $\af$ is global if and only if $A_g=A_{t(g)}$, for all $g\in \G$. Global actions of $\G$ on a ring $B$ induce by restriction partial actions of $\G$ on any ideal $A$ of $B$. Indeed, let $\beta=(B_g,\bt_g)_{g\in \G}$ be a global action of a groupoid $\G$ on a ring $B$, and $A\subseteq B$ an ideal of $B$. For each $x\in \G_0$ and $g\in \G$, consider $A_x=A\cap B_x$, $A_g=A_{t(g)}\cap \beta_g(A_{s(g)})$ and $\af_g=\beta_g|_{A_{g^{-1}}}$. Then $\af=(A_g,\af_g)_{g\in \G}$ is a partial action of $\G$ on $A$. Partial actions constructed in this way are called globalizable. The formal definition is given below.

Let $\af=(A_g,\alpha_g)_{g\in \G}$ be a partial action of a groupoid $\G$ on a ring $A$. A {\it globalization} of $\af$ is a pair $(B,\beta)$, where $B$ is a ring and $\beta=(B_g,\beta_g)_{g\in \G}$ is a global action of $\G$ on $B$, that satisfies
\begin{enumerate}[\hspace{.35cm}(G1)]
	\item  $A_x$ is an ideal of $B_x$, for each $x\in \G_0$,\smallbreak
	\item  $A_g=A_{t(g)}\cap \beta_g(A_{s(g)})$, for all $g\in \G$,\smallbreak
	\item  $\beta_g(a)=\alpha_g(a)$, for all $a\in A_{g^{-1}}$, \smallbreak
	\item  $B_g=\sum_{t(h)=t(g)}\beta_h(A_{s(h)})$, for all $g\in \G$.
\end{enumerate}

If $\alpha$ admits a globalization we say that $\alpha$ is {\it globalizable}. The globalization of $\alpha$ is unique, up to isomorphism; see section 2 of \cite{BP} for more details.

We recall that a partial action $\alpha=(A_g,\alpha_g)$ of $\G$ on $A$ is {\it unital} if every $A_g$ is unital, for $g\in \G$. In this case, the identity element of $A_g$,  denoted by $1_g$, is a central idempotent of $A$. By Theorem 2.1 of \cite{BP}, any unital partial action of a groupoid on a ring is globalizable.

For the convenience of the reader, we recall two useful properties of partial groupoid actions that were proved in Lemma 1.1 of \cite{BP}.

\begin{lem1}\label{two-proper}
	If $\alpha=(A_g,\alpha_g)$ is a partial action of a groupoid $\G$ on a ring $A$ then
	\begin{enumerate}[\hspace{0.2cm}\rm i)]
		\item $\alpha\m _g=\alpha_{g\m}$, for every $g\in \G$, and
		\item $\af_g(A_{g\m}\cap A_h)=A_g\cap A_{gh}$, for every $(g,h)\in \G_s\!\times\!\!\,_t\G $.
	\end{enumerate}
\end{lem1}

We end this section by noting that the study of partial actions of groupoids on rings can be reduced to the case of connected groupoids.

\begin{remark}{\rm Let $\G$ be a groupoid. Using the decomposition of $\G=\dot\cup_{X\in \G_0/\!\sim}\G_X$ in connected components, it is straightforward to check that partial actions of $\G$ on  a ring $\A$ induce by restriction partial actions of $\G_X$ on $\A$, for all $X\in \G_0/\!\!\sim$. Conversely, it is easy to see that partial actions of $\G$ on $A$ are uniquely determined by partial actions of $\G_X$, $X\in \G_0/\!\!\sim$, on $A$.}
\end{remark}

\section{Restriction and Extension Functors}

Throughout this section, $\G$ is a connected groupoid. Suppose that $x$ is an object of $\G$ and $\af=(A_g,\alpha_g)_{g\in \G}$ is a partial action of $\G$ on a ring $A$.
Notice that $\af$ induces by restriction a partial action
$(\A_h,\af_h)_{h\in\G(x)}$ of the group $\G(x)$ on the ring $\A_x$. Our main purpose in this section is to study the converse problem. More precisely, we will investigate when a partial action of $\G(x)$ on an ideal of $A$ can be extended to a partial action of $\G$ on $A$ .

\subsection{The categories $\Par$ and $\textsf{D}_\G(A)$}

Denote by $\Par$ the category  whose objects are  partial actions of $\G$ on  a fixed ring $\A$ and whose morphisms are defined as follows. Given  $\af=(\A_g,\af_g)_{g\in \G}$  and $\af'=(\A'_g,\af'_g)_{g\in \G}$ in  $\Par$, a morphism $\psi\colon \af\to \af'$ is a set of ring homomorphisms
$\psi=\{\psi_{y}:A_{y}\to A'_{y}\}_{y\in \G_0}$ such that, for every $g\in \G$,
\[\psi_{t(g)}( \A_g)\subseteq \A'_g,\qquad \af'_g\psi_{s(g)}(a)=\psi_{t(g)}\af_g(a),\,\,\text{ for all }\, a\in  \A_{g\m}.\]

In order to define the category $\textsf{D}_\G(A)$, we need some extra notation. Given $x\in \G_0$, consider on $\G(x,\,\,):=\{g\in \G\,:\,s(g)=x\}$ the following equivalence relation:
\begin{equation}\label{equiv} g\sim_x l\text{  if and only if  } t(g)=t(l).\end{equation}
 A transversal for $\sim_x$  such that $\tau_x=x$ will be called a \emph{transversal for $x$} and denoted by $\ta(x)$, that is,
$\tau(x)=\{\tau_{y}:y\in \G_0\}$ where $\ta_y$ is a chosen morphism in $\G(x,y)$, for each $y\in\G_0$, and we take $\ta_x=x$.

For a fixed transversal $\tau(x)$, we associate to each $\af=(\A_g,\af_g)_{g\in \G}\in\Par$ the triple $(A_{\alpha},\af_{\ta(x)},\alpha_{(x)})$, where
\begin{align}\label{3-tuple}
&\A_\af=\{\A_y\}_{y\in\G_0},& &\af_{\ta(x)}=\{\af_{\ta_{y}}\colon \A_{\ta\m_{y}}\to \A_{\ta_{y}}\}_{y\in \G_0},& &\af_{(x)}=(\A_h, \af_h)_{h\in \G(x)}.&
\end{align}
Observe that  $A_{\ta\m_y}$ is an ideal of $A_{t(\tau\m_y)}=A_x$ and $A_{\ta_{y}}$ is an ideal of $A_{t(\tau_{y})}=A_y$, for all $y\in \G_0$. Hence, $\alpha_{\tau_y}$ is a ring isomorphism from an ideal of $A_x$ to an ideal of $A_y$. Since $\tau_x=x$, it follows from (P2) that $\af_{\ta_{x}}=\af_x=\id_{A_x}$.
As in Proposition \ref{group:connec}, we consider $g_x:=\ta^{-1}_{t(g)}g\ta_{s(g)}$. Since $\tau_{t(g)}g_x\tau\m_{s(g)}=g$ and we want to recover $\af$ from the triple given in \eqref{3-tuple}, it is natural to investigate the codomain of the composition $\alpha_{\ta_{t(g)}}\af_{g_x}\af_{\tau\m_{s(g)}}$ in the inverse semigroup of partial bijections of $A$; see $\S\,$3.2 bellow to the definition of composition between partial bijections. By  Lemma \ref{two-proper}, this codomain is
\begin{align*}
\af_{\ta_{t(g)}}(A_{\ta\m_{t(g)}}\cap \af_{g_x}(A_{\ta\m_{s(g)}}\cap A_{g\m_x}))&=\af_{\ta_{t(g)}}(A_{\ta\m_{t(g)}}\cap A_{\ta\m_{t(g)}g}\cap A_{g_x})\\
&=A_{\ta_{t(g)}}\cap A_{g}\cap A_{g\ta_{s(g)}},
\end{align*}
which is an ideal of $A_{t(g)}$, for all $g\in \G$.

The triple that we introduced in $\eqref{3-tuple}$ suggests to consider the category $\Dat$ that we describe below. The objects of $\Dat$ are triples $(I,\gamma_{\ta(x)},\gamma_{(x)})$, where
\begin{enumerate}[$\quad \circ$]
	\item  $I=\{I_y\}_{y\in\G_0}$ is a family of ideals of $\A$, \smallbreak
		
	\item  $\gamma_{\ta(x)}=\{\gamma_{\ta_{y}}:\I_{\ta_{y}^{-1}}\to\I_{\ta_{y}}\}_{y\in\G_0}$ is a family of ring isomorphisms, $\I_{\ta_{y}^{-1}}$ is an ideal of $I_x$,
	 $\I_{\ta_{y}}$ is an ideal of $\I_y$ and $\gamma_{\ta_x}=\gamma_x=\id_{I_x}$,\smallbreak
	
	\item  $\gm_{(x)}=(\I_h,\gm_h)_{h\in\G(x)}$  is a  partial group action of $\G(x)$ on $\I_x$,
\end{enumerate}
that satisfy
\begin{align}\label{cond-ideal}
\gamma_{\ta_{t(g)}}(I_{\ta\m_{t(g)}}\cap \gamma_{g_x}(I_{\ta\m_{s(g)}}\cap I_{g\m_x}))\, \text{ is an ideal of }\, I_{t(g)}, \,\,\text{for all}\,\,g\in \G.
\end{align}
A morphism $f\colon (\I, \gamma_{\ta(x)}, \gm_{(x)})\to (\I', \gamma'_{\ta(x)}, \gm'_{(x)})$ in $\Dat$ is a set of ring homomorphisms $f=\{f_y:\I_y\to \I'_y\}_{y\in \G_0}$ such that
\begin{align*}
&f_y(\I_{\ta_{y}} )\subseteq \I'_{\ta_{y}}\,\,\text{ if }\,\  y\neq x,& &f_x(\I_{\ta\m_{y}} )\subseteq \I'_{\ta\m_{y}},& &\gamma'_{\ta_{y}}f_x=f_y\gamma_{\ta_{y}}.&
\end{align*}
Additionally we require that $f_x\colon \gm_{(x)} \to \gm'_{(x)}$ be a morphism of partial group actions, i.\,e., $f_x\colon \I_x\to \I'_x$ is a ring homomorphism that satisfies
\begin{align*}
&f_x(\I_h)\subseteq \I'_h\,\,\text{ and }\,\, \gm'_hf_x(a)=f_x\gm_h(a), \,\text{ for all }\, h\in \G(x), \,\,a\in  \I_{h\m}.&
\end{align*}

Now, we prove that the category $\Dat$ does not depend neither on the choice of the object $x\in\G_0$, nor on the choice of the transversal $\ta(x)$ of $x$.

\begin{prop1}\label{independent}
Let $x,z\in \G_0$,  $\ta(x)$ be a transversal for $x$ and $\,\lb(z)$ a transversal for $z$. Then $\Dat$ and $\NDat$ are isomorphic categories.
\end{prop1}

\begin{proof}
We start by defining the functor $F^{\lb(z)}_{\ta(x)}:\Dat\to \NDat$; the converse $F^{\ta(x)}_{\lb(z)}$ is defined similarly. For the object  $\gamma=(I,\gamma_{\ta(z)},\gamma_{(x)})\in\Dat$ we associate the object $\gamma'=(I',\gamma'_{\lambda(z)},\gamma'_{(z)})\in\NDat$ chosen in the following way:
	\begin{enumerate}[$\quad\circ $]
		\item  if $I=\{I_y\}_{y\in\G_0}$ we take $I'=\{I'_y\}_{y\in\G_0}$ such that $I'_z=I_x$, $I'_x=I_z$ and $I'_y=I_y$ for all $y\notin\{ x, z\}$,
		\item  if $\gamma_{\ta(x)}=\{ \gamma_{\ta_{y}}:I_{\ta\m_y}\to I_{\ta_y}\}_{y\in\G_0}$ we take $\gamma'_{\lb(z)}=\{ \gamma'_{\lb_{y}}:I'_{\lb\m_y}\to I'_{\lb_y}\}_{y\in\G_0}$
		such that $\gamma'_z=\gamma'_{\lb_z}=\gamma_{\ta_{(x)}}=\gamma_x$, $\gamma'_{\lb_z}=\gamma_{\ta_z}$ and $\gamma'_{\lb_y}=\gamma_{\ta_y}$
		for all $y\notin\{ z, x\}$,
		\item  if $\gamma_{(x)}=(I_h,\gamma_h)_{h\in\G(x)}$ we take  $\gamma'_{(z)}=(I'_l,\gamma'_l)_{l\in\G(z)}$ such that $I'_l=I_{\phi(l)}$, $\gamma'_l=\gamma_{\phi(l)}$ for
		all $l\in\G(z)$, where $\phi:\G(z)\to\G(x)$ is the group isomorphism defined by $l\mapsto \ta\m_zl\ta_z$.
	\end{enumerate}
Given a morphism $f=\{f_y\}_{y\in \G_0}:\gamma\to\delta$ in $\Dat$, we associate to it the morphism $f'=\{f'_y\}_{y\in \G_0}:\gamma'\to\delta'$ in $\NDat$, where $f'_{z}=f_{x}$, $f'_x=f_z$ and $f'_y=f_y$ for all $y\notin\{ x, z\}$.
This way we obtain the functors $F_{\ta(x)}^{\lb(z)}$ and $F_{\lb(z)}^{\ta(x)}$. It is straightforward to check that the compositions
	$F_{\lb(z)}^{\ta(x)}\circ F_{\ta(x)}^{\lb(z)}$ and $F_{\ta(x)}^{\lb(z)}\circ F_{\lb(z)}^{\ta(x)}$ are respectively the identity functors
	of $\Dat$ and $\NDat$.
\end{proof}

\noindent{\bf Notation:} Whenever $x\in \G_0$ and the transversal $\tau(x)$ for $x$ are fixed, we will denote the category $\Dat$ simply by $\textsf{D}_\G(A)$, that is, without any mention to the transversal $\tau{(x)}$.

\subsection{The functors $\Res$ and $\Ext$}
In the previous subsection we constructed the category $\DA(A)$ using the known information
about the objects of $\Par$. We will investigate if the information extracted from $\Par$ is enough to reverse the process. In order to do this, we will define functors that relate these categories.

Fix  $x\in \G_0$ and $\tau(x)$ a transversal for $x$. Let $\af=(\A_g,\af_g)_{g\in \G}$ be an object in  $\Par$ and $\psi:\af'\to \af''$ a morphism in $\Par$. The association
$\Res\colon \Par\to \DA(A)$ given by
\begin{align*}
&\Res(\af)=(\A_\af, \af_{\ta(x)}, \af_{(x)}),& &\Res(\psi)=\psi,&
\end{align*}
is a functor. Since $\Res(\alpha)$ is the restriction of the partial groupoid action of $\G$ on $A$ to the partial group action of $\G(x)$ on $A_x$, we will say that $\Res$ is the {\it restriction functor}.

In order to construct the reverse functor, we start by recalling the definition of composition of partial bijections. Let $Z$ be a nonempty set, $f:X\to Y$ and $f':X'\to Y'$ be partial bijections of $Z$. The composition $f'f$ is defined by $f'f:f^{-1}(X'\cap Y)\to f'(X'\cap Y)$ which is a partial bijection of $Z$.

Given an object $(\I, \gamma_{\ta(x)},\gm_{(x)})$ in $\DA(A)$ we set $\theta=(B_g,\theta_g)_{g\in \G}$, where each $\theta_g$ is the following ring isomorphism
\begin{equation}\label{afcon}\theta_g=\begin{cases}
\id_{I_{y}},& \text{ if } g=y\in \G_0,\\
\gm_{\ta_{t(g)}}\gamma_{g_x}\gamma\m_{\ta_{s(g)}},& \text{ if } g\notin \G_0,
\end{cases}
\end{equation}
and each $B_g$ is taken as the range of $\theta_g,$ for all $g\in \G$.  Notice that $B_y=I_y$ for all $y\in \G_0$ and
\begin{align}
\label{dom} B_{g}=
\gamma_{\ta_{t(g)}}(I_{\ta\m_{t(g)}}\cap \gamma_{g_x}(I_{\ta\m_{s(g)}}\cap I_{g\m_{x}})),& \text{ if } g\notin \G_0.
\end{align}
Hence $B_y$ is an ideal of $A$, for all $y\in \G_0$. Also, by \eqref{cond-ideal}, $B_g$ is an ideal of $B_{t(g)}=I_{t(g)}$, for all $g\notin \G_0$.

\begin{teo1}\label{outro} The pair $\theta=(B_{g},\theta_{g})_{g\in \G}$ constructed above  is an object of $\Par$.
\end{teo1}
\begin{proof}
Let $(g,h)\in \G_s\times\!_t\G$ such that  $g, h\notin \G_0$. Then
\begin{align*}
\theta_{g}\theta_{h}&=\gamma_{\ta_{t(g)}}\gamma_{g_x}\gamma\m_{\ta_{s(g)}}\gamma_{\ta_{t(h)}}\gamma_{h_x}\gamma\m_{\ta_{s(h)}}\\
&=\gamma_{\ta_{t(g)}}\gamma_{g_x}\gamma_{h_x}\gamma\m_{\ta_{s(h)}}
\end{align*}
Since $\gamma_{(x)}$ is a partial action of $\G(x)$ on $I_x$, we have that $\gamma_{(gh)_x}=\gamma_{g_xh_x}$ is an extension of $\gamma_{g_x}\gamma_{h_x}$, that is, $\gamma_{g_x}\gamma_{h_x}\leq \gamma_{g_xh_x}$. Thus
\begin{align*}
\theta_{g}\theta_{h}&\leq \gamma_{\ta_{t(g)}}\gamma_{g_xh_x}\gamma\m_{\ta_{s(h)}}\\
&=\gamma_{\ta_{t(gh)}}\gamma_{(gh)_x}\gamma\m_{\ta_{s(gh)}}\\
&=\theta_{gh}.\end{align*}
If $h=y\in \G_0$  and $s(g)=y$ then $\theta_{g}\theta_h=\theta_g=\theta_{gh}$. Similarly, $\theta_{g}\theta_h=\theta_{gh}$ when $g=y=t(h)\in \G_0$ or $g=h=y\in \G_0$. By $\S\,$3.1 of \cite{BPP}, the condition $\theta_g\theta_h\leq \theta_{gh}$ for all $(g,h)\in \G_s\times\!_t\G$ is equivalent to the conditions (P3) and (P4) given in the definition of partial groupoid action. Hence, the result follows.
\end{proof}

Let $\gamma=(\I, \gamma_{\ta(x)},\gm_{(x)})$ be an object in $\DA(A)$ and $f:\gamma'\to \gamma''$ a morphism in $\DA(A)$.
By Theorem \ref{outro}, $\Ext(\gamma):=\theta=(B_g,\theta_g)_{g\in\G}\in \Par$. Also, if $f=\{f_y:I'_y\to I''_y\}_{y\in \G_0}$ then we define $\Ext(f):=\psi:\theta'\to \theta''$ by $\psi_y=f_y$, for each $y\in \G_0$. Hence, the association $\Ext:\DA(A)\to\Par$ given by
\begin{align*}
&\Ext(\gamma)=\theta=(B_g,\theta_g)_{g\in\G},& &\Ext(f)=\psi,&
\end{align*}
is a functor. The functor $\Ext$ will be called the {\it extension functor}.

\vd

The next result establishes the relation between the restriction functor $\Res$  and  the extension functor $\Ext$.

\vu

\begin{prop1}\label{sub} Let $\af=(A_g,\af_g)_{g\in\G}\in  \Par$ and $\theta:=\Ext(\Res(\af))$. Then
\begin{enumerate}[\hspace{0.2cm}\rm i)]
\item $\theta\leq \af$, that is, $\af_g$ is an extension of $\theta_{g}$, for all $g\in \G$, \smallbreak

\item $\theta=\af$ if and only if $A_g\subseteq A_{\ta_{t(g)}}\cap A_{g\ta_{s(g)}}$, for all $g\notin\G_0$, \smallbreak

\item  $\Res(\Ext(\gamma))=\gamma$, for all $\gamma\in\DA(A)$.
\end{enumerate}
\end{prop1}

\begin{proof} Observe that for $g\notin \G_0$ we have
\begin{align*}
\,\,\,\theta_g &\overset{\mathclap{\eqref{afcon}}}{=}\af_{\ta_{t(g)}}\af_{g_x}\af_{\ta\m_{s(g)}}&\\
      &=\af_{\ta_{t(g)}}\af_{\ta\m_{t(g)}g\ta_{s(g)}}\af_{\ta\m_{s(g)}}& \\
      &\overset{\mathclap{(\ast)}}{\leq}\af_{\ta_{t(g)}\ta\m_{t(g)}g\ta_{s(g)}\ta\m_{s(g)}}&\\
      &=\af_g,
\end{align*}
where $(\ast)$ follows from Remark \ref{obs-inicial}.
If $g=y\in \G_0$, then $\theta_g=\theta_y=\af_y=\af_g$ and i) follows. For ii), notice that
\begin{align*}
B_g &\overset{\mathclap{\eqref{dom}}}{=}\af_{\ta_{t(g)}}(A_{\ta\m_{t(g)}}\cap\af_{g_x}(A_{\ta\m_{s(g)}}\cap A_{g\m_x}))\\
&=\af_{\ta_{t(g)}}(A_{\ta\m_{t(g)}}\cap A_{g_x\ta\m_{s(g)}}\cap A_{g_x})\qquad\qquad (\text{by Lemma \ref{two-proper}})\\
&=\af_{\ta_{t(g)}}(A_{\ta\m_{t(g)}}\cap A_{\ta\m_{t(g)}g}\cap A_{g_x})\\
&= A_{\ta_{t(g)}}\cap A_g\cap A_{\ta_{t(g)}g_x} \qquad\qquad (\text{by Lemma \ref{two-proper}})\\
&=A_g\cap(A_{\ta_{t(g)}}\cap A_{g\ta_{s(g)}}),
\end{align*}
for all $g\notin \G_0$, and whence we obtain ii). The item iii) is easily checked using the definitions of $\Ext$ and $\Res$.
\end{proof}

\begin{cor1} The functors $\Ext$ and $\Res$ form an adjoint pair.
\end{cor1}
\begin{proof} Denote the functors $\Ext$ and $\Res$ respectively by $\texttt{E}$ and $\texttt{R}$. By Proposition \ref{sub} (iii),  $\eta:{\texttt{Id}} _{\DA(A)}\to \nres\circ\next$, given by $\eta_{\gamma}={\rm id}_\gamma$, $\gamma\in \textsf{D}_\G(A)$, is a natural transformation. Let $\af=(A_g,\af_g)_{g\in \G}\in \Par$ and 
	$\theta_\alpha=\next(\nres(\af))=(B_g,(\theta_\alpha)_g)_{g\in \G}$. By Proposition \ref{sub} (i), $B_y\subset A_y$ for all $y\in \G_0$. Then,
	$\varepsilon\colon \next\circ \nres\to {\texttt{Id}}_{\Par} $ given by $\varepsilon_\alpha:\theta_\alpha\to \alpha$, is a natural transformation, where 
	$\varepsilon_\alpha=\{\varepsilon_{\alpha,y}:B_y\to A_y\}_{y\in \G_0}$ is the family of inclusions.
	
	Let $\gamma\in \DA(A)$. Note that $\theta_{\next(\gamma)}=\next(\nres(\next(\gamma)))=\next(\gamma)$ and whence  $\varepsilon_{\next(\gamma)}:\next(\gamma)\to\next(\gamma)$ is the identity morphism ${\rm id}_{\next(\gamma)}$. Consequently $\varepsilon_{\next(\gamma)}\circ \next(\eta_\gamma)=\varepsilon_{\next(\gamma)}\circ {\rm id}_{\next(\gamma)}={\rm id}_{\next(\gamma)}$. Since $\nres(\theta_\alpha)=\nres(\alpha)$, it is easy to check that $\nres(\varepsilon_\alpha)\circ \eta_{\nres(\alpha)}={\rm id}_{\nres(\alpha)}$ for all $\alpha\in \Par$. It follows from Theorem 2 (p.83) of \cite{ML} that $(\next,\nres )$ is an adjoint pair.
\end{proof}

Now we prove a result that will be useful in the rest of this subsection.

\begin{prop1}\label{pro:equivalentes} Let $\af=(A_g,\af_g)_{g\in\G}\in  \Par$. Then the following conditions are equivalent:
	\begin{enumerate} [\hspace{0.2cm}\rm i)]
		\item $A_g\subseteq A_{\ta_{t(g)}}\cap A_{g\ta_{s(g)}}$, for all $g\in \G\setminus\G_0$,\smallbreak
		\item $A_g\subseteq A_{\ta_{t(g)}}$, for all $g\in \G\setminus\G_0$.
	\end{enumerate}
\end{prop1}
\begin{proof}
	
	Clearly  i) $\Rightarrow$ ii). Conversely, if $g\in \G\setminus\G_0$ then  $A_{g\m}\subseteq A_{\ta_{t(g\m)}}=A_{\ta_{s(g)}}$. Thus, $$\af_{g\m}(A_{g})=A_{g\m}=A_{g\m}\cap A_{\ta_{s(g)}}=\af_{g\m}(A_{g}\cap A_{g\ta_{s(g)}}).$$ Since $\af_{g\m}$ is a ring isomorphism, we obtain that
	$A_g=A_g\cap A_{g\ta_{s(g)}}$ and consequently $A_g\subseteq A_{g\ta_{s(g)}}$.
\end{proof}

Suppose that $\af=(A_g,\af_g)_{g\in\G}\in \Par$ satisfies one of the equivalent conditions of Proposition \ref{pro:equivalentes}. Then, by Proposition \ref{sub}, $\af=\Ext(\Res(\af))$.
Motivated by this, we introduce the following definition. \smallbreak

A partial action $\af=(A_g,\af_g)_{g\in\G}\in \Par$ will be called {\it recoverable by $\Ext$} if there is an object $x\in \G_0$ and a transversal $\ta(x)=\{\tau_y:\,y\in \G_0\}$ for $x$ such that
\begin{align}\label{recover}
&A_g\subset A_{\tau_{t(g)}},\,\,\,\,\text{for all}\,\,g\in \G\setminus\G_0.&
\end{align}

In general not all partial actions are recoverable by $\Ext$; see Example \ref{brandt2}. However, we will see bellow that an important class of partial groupoid actions are recoverable by $\Ext$.\smallbreak

We recall from \cite{BPP} that a partial action $\af=(A_g,\af_g)_{g\in \G}$ of $\G$ on a ring $A$ is {\it  group-type} if there exist $x\in \G_0$ and a transversal $\tau(x)=\{\tau_{y}:y\in \G_0\}$ for $x$ such that
\begin{align}
\label{cond1} A_{\tau\m_y}=A_x \ \ \text{and} \ \ A_{\tau_y}=A_y, \ \ \text{ for all } \ y\in\G_0.
\end{align}
 Denote by  $\GT$ the full subcategory of $\Par$  whose objects are the group-type actions of $\G$ on $A$. Global actions are examples of group-type actions.

 Also, if $x\in \G_0$ and $\ta(x)$ is a fixed transversal for $x$ then the full subcategory of $\DA(A)$ with objects  $(\I, \gamma_{\ta(x)},\gm_{(x)})$ satisfying $I_{\ta_y}=I_{y}$ and $I_{\ta\m_y}=I_{x}$, for all $y\in \G_0$, will be denoted by $\GDA(A)$.

\vu
\begin{cor1} The following statements hold:
	\begin{enumerate}[\hspace{0.2cm}\rm i)]
		\item if $\af\in \GT$ then $\af$ is recoverable by $\Ext$,\smallbreak
		\item the categories $\GT$ and $\GDA(A)$ are isomorphic.
	\end{enumerate}
\end{cor1}
\begin{proof}
Consider $\af=(A_g,\af_g)_{g\in \G}\in \GT$. Then there exist $x\in \G_0$ and $\ta(x)$ a transversal for $x$ such that \eqref{cond1} holds. Hence  $A_g\subset A_{t(g)}=A_{\ta_{t(g)}}$, for all $g\in \G\setminus \G_0$. Thus $\af$ is recoverable by $\Ext$ and i) follows. Notice that Proposition \ref{pro:equivalentes} ii) and  Proposition \ref{sub} ii) imply that $\Ext(\Res(\af))=\af$. Then, from Proposition \ref{sub} iii) we conclude that $\GT$ and $\GDA(A)$ are isomorphic.
\end{proof}

\begin{remark}{\rm  The objects of $\GT$ are an important subclass of partial groupoid actions, as can be seen in \cite{BPP}. Particularly, it was proved in section 4 of \cite{BPP} that if $\af\in \GT$ then the partial skew groupoid ring $A\star_\af\G$ can be realized as a partial skew group ring. In this sense, some problems involving partial groupoid actions are reduced to analogous problems involving partial group actions.}
	
\end{remark}

\subsection{Examples}
We start this subsection with some examples that illustrate the construction of a partial groupoid action from an object of $\DA(A)$.

\begin{exe}\label{hexd} {\rm Let $\G=\{g,h,l,m,l^{-1},m^{-1}\}$ be the groupoid with objects $\G_0=\{x,y\}$ and the following composition rules
		\[ g^2=x,\quad h^{2}=y,\quad lg=m=hl,\quad g\in\G(x),\quad h\in\G(y)\,\,\text{ and }\,\, l,m\in\G(x,y). \]
		The diagram bellow illustrates the structure of $\G$:
		\[\xymatrix{& x\ar[r]^{l}  &y \ar[d]^{h}\\
			& x\ar[u]^{g} \ar[r]^{m} & y} \]
		Take the transversal $\ta(x)=\{\ta_x=x, \ta_y=l\}$ for $x$. Let $A$ be a ring and $J,L,K$ ideals of $A$ such that $L\subset K$ and $J\subset K$. Consider two ring automorphisms $\gamma$ and $\sigma$ of $K$ such that $\sigma(L)=L$, $\sigma|_L=\sigma\m|_L$ and $\gamma(J)=J$.  Then we take the following object $\big(I,\gamma_{\ta(x)},\gamma_{(x)}\big)$ in $\DA(A)$:
		
		\begin{enumerate}[$\quad\circ $]
			\item $I=\{I_x,\,I_y\}$ with $I_x=I_y=K$,\smallbreak
			
			\item  $\gamma_{\ta(x)}=\{\gamma_{\ta_x}, \gamma_{\ta_y}\}$ with  $\gamma_{\ta_x}=\gamma_x=\id_{I_x}=\id_K$ and  $\gamma_{\ta_y} =\gamma_l=\gamma$, \smallbreak
			
			\item  $\gamma_{(x)}$ is the partial action of $\G(x)=\{x, g\}$ on $I_x=K$ given by $\gamma_x=\id_{K}$ and $\gamma_g=\gamma_{g\m}=\sigma:L\to L$.
		\end{enumerate}	
		Then $\theta=\Ext\big(I,\gamma_{\ta(x)},\gamma_{(x)}\big)$ is the following partial action of $\G$ on $A$. First, note that
		$$g_x=h_x=m_x=m\m_x=g\,\,\,\text{ and }\,\,\, l_x=l\m_x=x_x=y_x=x.$$
		From \eqref{afcon} and \eqref{dom} it follows that
		\begin{align*}
		&\theta_x=\theta_y=\id_K,\qquad \theta_g=\sigma|_{L}=\theta_{g\m},\qquad \theta_h=\gamma\sigma\gamma\m|_{\gamma(J\cap\sigma(J\cap L))}=\theta_{h\m},\\[.5em]
		& \theta_m=\gamma\sigma:\sigma(J\cap L)\to \gamma(J\cap L), \quad\theta_{m\m}=\sigma\gamma\m:\gamma(J\cap L)\to \sigma(J\cap L).
		\end{align*}}		
\end{exe}

\begin{exe} \label{brandt} {\rm Let $m$ be a positive integer and $G$ a group. Denote by $\Gamma_G^m$ the following connected groupoid associated to $m$ and $G$. The set of
		objects of $\Gamma_G^m$ is  $\mathbb{I}_m=\{1,\ldots, m\}$. A morphisms in $\Gamma_G^m$ is an element $(i,g,j)$, where $g\in G$ and $i,j\in\mathbb{I}_m$. Also, $s(i,g,j)=j$, $t(i,g,j)=i$  and the composition  is given by  \[(i,g,j)(j,h,l)=(i,gh,l),\,\text{ for all } i,j,l\in\mathbb{I}_m\,\text{ and }\, g,h\in G.\]
		Fix $i_0\in \mathbb{I}_m$ and the transversal $\tau(i_0)=\{\tau_i=(i,e,i_0)\}_{i\in \mathbb{I}_m}$ for $i_0$, where $e$ denotes the identity element of $G$. Notice that $\Gamma_G^m(i_0)=\{(i_0,g,i_0)\mid g\in G\}\simeq G$. Let
		$A$ be a unital ring and $\{e_i\}_{i\in \mathbb{I}_m}$ a set of orthogonal idempotents in the center of $A$ whose sum is $1_A$.  Suppose that $G$ acts on $Ae_{i_0}$ by ring isomorphisms $\gamma_g$, $g\in G$, and there exists a family of ring isomorphism $\gamma_{\ta_i}:A{e_{i_0}}\to A{e_{i}}$, for each $i\in \mathbb{I}$. Then, we take $\big(I,\gamma_{\ta(i_0)},\gamma_{(i_0)}\big)\in \DA(A)$ in the following way: }
	
	\vu
	
	\begin{enumerate}[$\quad\circ $]
		\item {\rm $I=\{I_i\}_{i\in\mathbb{I}_m}$ where $I_i=Ae_{i},$} \smallbreak

		\item {\rm $\gamma_{\ta({i_0})}=\{\gamma_{\ta_i}\}_{i\in\mathbb{I}_m}$,} \smallbreak
		
		\item  {\rm $\gamma_{(i_0)}=\{\gamma_{g}\}_{g\in G}$, where $\gamma_g$ is the action of $G$ on $A{e_{i_0}}$.}
	\end{enumerate}	
	{\rm Then $\theta=\Ext\big(I,\gamma_{\ta(i_0)},\gamma_{(i_0)}\big)$ is the following global action of $\Gamma_G^m$ on $A$:}
	\[\theta_{(i,g,i)}=\id_{Ae_i},\qquad \theta_{(i,g,j)}=\gamma_{\ta_{i}}\gamma_{g}\gamma\m_{\ta_{j}}:Ae_j\to Ae_i, \qquad i,j\in\mathbb{I}_m,\,\,g\in G.  \]	
\end{exe}

\begin{exe}\label{ex:datum-global}
	{\rm Let $\G=\{g,h,l,m,l^{-1},m^{-1}\}$ be the groupoid given in Example \ref{hexd} and  $\ta(x)=\{\ta_x=x, \ta_y=l\}$. Let $A$ be a ring, $\sigma,\gamma$ ring automorphisms of $A$ with $\sigma\m=\sigma$ and $e\in A$ a central idempotent of $A$. Consider $\big(I,\gamma_{\ta(x)},\gamma_{(x)}\big)\in \DA(A)$, where
		\begin{enumerate}[$\quad\circ $]
			\item  $I=\{I_x,\,I_y\}$ with $I_x=Ae$  and $I_y=A\gamma(e)$,\smallbreak
			\item $\gamma_{\ta(x)}=\{\gamma_{\ta_x}, \gamma_{\ta_y}\}$ with  $\gamma_{\ta_x}=\gamma_x=\id_{I_x}$ and $\gamma_{\ta_y}=\gamma_l=\gamma|_{I_x}:I_x\to I_y$,  \smallbreak
			\item  $\gamma_{(x)}$ is the partial action of $\G(x)=\{x, g\}$ on $I_x$ given by $\gamma_x=\id_{I_x}$ and \linebreak $\gamma_g=\gamma_{g\m}=\sigma|_{Ae\sigma(e)}:Ae\sigma(e)\to Ae\sigma(e)$.
		\end{enumerate}	
		Then $\theta=\Ext\big(I,\gamma_{\ta(x)},\gamma_{(x)}\big)$ is the partial action of $\G$ on $A$ given by
		\begin{align*}
		&\theta_x=\id_{I_x},\qquad \theta_y=\id_{I_y},\qquad \theta_g=\gamma_g,\qquad \theta_l=\gamma_l,\qquad\theta_m=\gamma\sigma|_{Ae\sigma(e)},\\[.2em]
		&\theta_h=\gamma\sigma \gamma\m|_{A\gamma(e\sigma(e))},\qquad\theta_{l\m}=\gamma_l\m,\qquad\theta_{m\m}=\sigma\gamma\m|_{A\gamma(e\sigma(e))}.
		\end{align*}}
\end{exe}

\smallbreak

We finish this subsection with an example of a partial groupoid action that is not recoverable by $\Ext$.

\begin{exe} \label{brandt2}
	{ \rm  Consider the groupoid $\G$ given in Example \ref{hexd} and $A=\sum_{i=1}^6\Bbbk e_i$, where $\{e_i\}_{i=1}^6$ is a set of orthogonal idempotents whose sum is $1_A$ and $\Bbbk$ is a field. Then we have the following partial action $\af=(A_g,\af_g)_{g\in \G}$ of  $\G$ in $A$:
		\begin{align*}
		&\af_{x}=\id_{A_x},\,\,\, \af_{y}=\id_{A_y},\,\, \text{ where}\,\, A_{x}=\Bbbk e_1\oplus\Bbbk e_2\oplus\Bbbk e_3,\,\,A_{y}=\Bbbk e_4\oplus\Bbbk e_5\oplus\Bbbk e_6,&\\[.2em]
		&\af_{g}=\id_{A_g},\,\,\, \af_{h}=\id_{A_h},\,\, \text{ where}\,\,A_{g}=\Bbbk e_1,\quad A_h=\Bbbk e_6,&\\[.2em]
		&\af_l:\Bbbk e_2\to \Bbbk e_4,\,\,\,\af_l(\lambda e_2)=\lambda e_4\quad \af_m:\Bbbk e_3\to \Bbbk e_5,\,\,\,\af_m(\lambda e_3)=\lambda e_5,\quad \lambda\in \Bbbk.&
		\end{align*}	
It is easy to check that $\af$ is a partial action of $\G$ on $A$. If $\ta(x)=\{\tau_x=x,\ta_y=l\}$ then notice that $A_h\nsubseteq A_{\ta_{t(h)}}=A_l$. The other possibility of a transversal for $x$ is $\ta(x)=\{\tau_x=x,\ta_y=m\}$. In this case, we have again that $A_h\nsubseteq A_{\ta_{t(h)}}=A_m$. In both cases, \eqref{recover} is not satisfied. Similarly, the condition \eqref{recover} is not satisfied for any transversal $\tau(y)$. Hence, $\af$ is not recoverable by $\Ext$.}
\end{exe}

\section{Globalization}

The problem of investigating when a partial action is globalizable, that is, when it can be obtained as restriction of a global action, is relevant
for partial actions; see \cite{BP} (resp. \cite{DE} ) for details on globalization of partial groupoid (resp. group) actions.

\vu

Throughout this section, $\G$ denotes a connected groupoid, $x\in \G_0$ is an object of $\G$ and  $\tau(x)=\{\tau_{y}:y\in \G_0\}$ is a transversal for $x$. We fix a ring $A$, $\gamma=(\I, \gamma_{\ta(x)},\gm_{(x)})\in \DA(A)$ and $\theta=\Ext(\gamma)$ the partial action of $\G$ on $A$ given by \eqref{afcon} and \eqref{dom}. \vu

Our main purposes here are the following. First, we will relate the existence of the globalization for $\theta$ with the existence of the globalization for $\gamma_{(x)}$. Second, we will investigate how to construct a globalization for $\theta$ from a globalization of $\gm_{(x)}$.

\vu

\begin{prop1}\label{prop:globa} Assume that $I_y$ is a unital ring, for all $y\in \G_0$. Then the  partial groupoid action $\theta$ of $\G$ on $A$ is globalizable if and only if the partial group action $\gm_{(x)}$ of $\G(x)$ on $I_x$ is globalizable and $I_{\ta\m_{y}}$ is a unital ring, for all $y\in \G_0$.
\end{prop1}
\begin{proof} Suppose that $\theta$ is globalizable. By Theorem 2.1 of \cite{BP}, $B_g$ is a unital ring, for all $g\in \G$. Let $1_g$ be a central idempotent of $A$ such that $B_g=A1_g$, $g\in \G$.	
Then $I_h=B_h=A1_h=A1_x1_h=B_x1_h=I_x1_h$, for all $h\in \G(x)$. Thus, by Theorem 4.5 of \cite{DE}, the  partial group action $\gm_{(x)}$ of $\G(x)$ on $I_x$ is globalizable. Also, $I_{\ta\m_{y}}=B_{\ta\m_{y}}=A1_{\ta\m_{y}}$, for all $y\in \G_0$.
	
	Conversely, since $\gm_{(x)}$ is globalizable it follows from Theorem 4.5 of \cite{DE} that $I_h$ is a unital ring, for all $h\in \G(x)$. Hence there exist central idempotents $1_h$ of $I_x$ such that $I_h=I_x1_h$. Consider also central idempotents $1_y$ of $A$ and $1_{\ta\m_{y}}$ of $I_x$ such that $I_y=A1_y$ and $I_{\ta\m_{y}}=I_x1_{\ta\m_{y}}$ for all $y\in \G_0$. Notice that $a1_{\ta\m_y}=a(1_x1_{\ta\m_y})=(a1_x)1_{\ta\m_y}=1_{\ta\m_y}(a1_x)=1_{\ta\m_y}(1_xa)=1_{\ta\m_y}a$, for all $a\in A$. Thus $1_{\ta\m_y}$ is a central idempotent of $A$. Then $1_g=\gamma_{\ta_{t(g)}}(1_{\ta\m_{t(g)}}\gamma_{g_x}(1_{\ta\m_{s(g)}}1_{g\m_x}))$ is a central idempotent of $A$, for all $g\in \G\setminus \G_0$. It is straightforward to check that 
	$B_{g}=\gamma_{\ta_{t(g)}}(\I_{\ta\m_{t(g)}}\cap \gamma_{g_x}(\I_{\ta\m_{s(g)}}\cap \I_{g\m_{x}}))$ satisfies $B_g=A1_g$.
	 If $g=y\in \G_0$ then $B_g=B_y=I_y=A1_y$. Hence, Theorem 2.1 of \cite{BP} implies that $\theta$ is globalizable.	
\end{proof}
In general, it is not easy to calculate the globalization of a partial groupoid action (see $\S\,2$ of \cite{BP}). However, we can do it if $\theta=\Ext(\gamma)$ and the object $\gamma\in \DA(A)$ satisfies the following conditions:
\begin{enumerate}[$\quad$(C1)]
	\item $\gamma\in \GDA(A)$, that is,  $I_{\ta\m_{y}}=I_x$ and $I_{\ta_{y}}=I_y$ for all $y\in \G_0$, \smallbreak
	\item $I_x$ is a unital ring and the partial group action $\gm_{(x)}$ of $\G(x)$ on $I_x$ admits a globalization $(J_x,\tilde{\gm}_{(x)})$, \smallbreak
	\item there are a ring $B$, a family of ideals $\{J_y\}_{y\in\G_0}$ of $B$ and a family of ring isomorphisms $\{\tilde{\gamma}_{\ta_y}:J_x\to J_y\}_{y\in\G_0}$ such that $I_y$ is an ideal of $J_y$ and $\tilde{\gamma}_{\ta_y}|_{I_{x}}=\gamma_{\ta_{y}}$, for all $y\in \G_0$.
\end{enumerate}
\vu

Suppose that $\gamma$ satisfies (C1) and (C2). Then $\gamma_{\ta_y}: I_x\to I_y$ is an isomorphism and (C2) implies that $I_y$ is unital. Thus $I_x=I_{\ta\m_y}$ and $I_y=I_{\ta_y}$ are unital rings, for all $y\in\G_0$. Since $\gamma_{(x)}$ is globalizable, it follows from Proposition \ref{prop:globa} that $\theta$ is globalizable. The next Proposition describes the globalization of $\theta$ if $\gamma$ also satisfies  (C3).

\begin{prop1}\label{teo:globaliza1} If $\gamma \in \DA(A)$ satisfies {\rm (C1)}, {\rm (C2)} and {\rm (C3)} above then the groupoid action $\beta$ of $\G$ on $B$ given by
	\[\bt=(C_g,{\bt}_g)_{g\in	\G},\qquad {C}_g=J_{t(g)},\qquad {\bt}_g=\tilde{\gamma}_{\ta_{t(g)}}\tilde{\gamma}_{g_x}\tilde{\gamma}\m_{\ta_{s(g)}}\,:\,J_{s(g)}\to J_{t(g)},\]
	is the globalization of $\theta$.
\end{prop1}
\begin{proof}
	Clearly ${\beta}$ is a global action of $\G$ on $B$ and ${\beta}_g(a)=\theta_g(a)$, for all $a\in B_{g^{-1}}$ and $g\in \G$. Notice that
	\begin{align*}
	B_{t(g)}\cap {\bt}_g(B_{s(g)})&= I_{t(g)}\cap{\bt}_g(I_{s(g)})\\
	&= I_{t(g)}\cap\tilde{\gamma}_{\ta_{t(g)}}\tilde{\gamma}_{g_x}\tilde{\gamma}\m_{\ta_{s(g)}}(I_{s(g)})\\
	&= I_{\ta_{t(g)}}\cap\tilde{\gamma}_{\ta_{t(g)}}\tilde{\gamma}_{g_x}\tilde{\gamma}\m_{\ta_{s(g)}}(I_{\ta_{s(g)}})\\
	&=\tilde{\gamma}_{\ta_{t(g)}}(I_x)\cap \tilde{\gamma}_{\ta_{t(g)}}\tilde{\gamma}_{g_x}(I_x)\\
	&=\tilde{\gamma}_{\ta_{t(g)}}(I_x\cap\tilde{\gamma}_{g_x}(I_x))\\
	&\overset{\mathclap{(\ast)}}{=}\tilde{\gamma}_{\ta_{t(g)}}(I_{g_x}),
	\end{align*}
	where $(\ast)$ holds because $(J_x,\tilde{\gm}_{(x)})$ is a globalization of the partial group action $\gm_{(x)}$ of $\G(x)$ on $I_x$. On the other hand,
	\begin{align*}
	B_g&=\gamma_{\ta_{t(g)}}(I_{\ta\m_{t(g)}}\cap\gamma_{g_x}(I_{\ta\m_{s(g)}}\cap I_{g\m_x}))\\
	&=\gamma_{\ta_{t(g)}}(I_{x}\cap\gamma_{g_x}(I_{x}\cap I_{g\m_x}))\\
	&=\gamma_{\ta_{t(g)}}(I_{g_x})\\
	&=\tilde{\gamma}_{\ta_{t(g)}}(I_{g_x}).
	\end{align*}
	Hence, $B_g=B_{t(g)}\cap {\bt}_g(B_{s(g)})$. Note also that, for each $g\in \G$, we have
	\begin{align*}
	\sum_{t(h)=t(g)}{\bt}_h(B_{s(h)})&=\sum_{t(h)=t(g)}\tilde{\gamma}_{\ta_{t(h)}}\tilde{\gamma}_{h_x}\tilde{\gamma}\m_{\ta_{s(h)}}(I_{s(h)})\\
	&=\tilde{\gamma}_{\ta_{t(g)}}\left(\sum_{t(h)=t(g)}\tilde{\gamma}_{h_x}(I_x)\right)\\
	&\overset{\mathclap{(\ast\ast)}}{=}\tilde{\gamma}_{\ta_{t(g)}}(J_x)\\
	&=J_{t(g)}\\
	&=C_g.
	\end{align*}
	In order to justify $(\ast\ast)$ note that, if $l\in \G(x)$ then $(\ta_{t(g)}l)_x=(\ta_{t(g)})_xl_x=xl=l$. Thus, $\sum_{t(h)=t(g)}\tilde{\gamma}_{h_x}(I_x)=\sum_{l\in \G(x)}\tilde{\gamma}_{l}(I_x)$. Since $(J_x,\tilde{\gm}_{(x)})$ is the globalization of $\gamma_{(x)}$ it follows that $J_x=\sum_{l\in \G(x)}\tilde{\gamma}_{l}(I_x)$.
\end{proof}

\vu

The next example illustrates the construction given by Proposition \ref{teo:globaliza1}.

\vu

\begin{exe}\label{ex:datum-globa}
	{\rm Let $\G$ be the groupoid given in Example \ref{hexd} and $\ta(x)=\{\ta_x=x, \ta_y=l\}$. Let $A$ be a ring, $\sigma,\gamma$ ring automorphisms of $A$ with $\sigma\m=\sigma$ and $e$ a central idempotent of $A$. Consider $\big(I,\gamma_{\ta(x)},\gamma_{(x)}\big)\in \DA(A)$ given by:
	\begin{enumerate}[$\quad\circ $]
		\item  $I=\{I_x,I_y\}$ with $I_x=Ae$ and $I_y=A\gamma(e)$,\smallbreak
		\item $\gamma_{\ta(x)}=\{\gamma_{\ta_x},\gamma_{\ta_y}\}$ with $\gamma_{\ta_x}=\gamma_x=\id_{I_x}$ and $\gamma_{\ta_y}=\gamma_l=\gamma|_{I_x}:I_x\to I_y$, \smallbreak
		\item  $\gamma_{(x)}$ is the partial action of $\G(x)=\{x, g\}$ on $I_x$ given by $\gamma_x=\id_{I_x}$ and \linebreak $\gamma_g=\gamma_{g\m}=\sigma|_{Ae\sigma(e)}:Ae\sigma(e)\to Ae\sigma(e)$.
	\end{enumerate}	
Then $\theta=(B_z,\theta_z\big)_{z\in\G}=\Ext(\gamma)$ is the partial action of $\G$ on $A$ given by
\begin{align*}
&\theta_x=\id_{I_x},\quad \theta_y=\id_{I_y},\quad \theta_g=\gamma_g,\quad \theta_l=\gamma_l,\quad\theta_m=\gamma\sigma|_{Ae\sigma(e)}:Ae\sigma(e)\to A\gamma(e\sigma(e)),\\
&\theta_h=\gamma\sigma \gamma\m|_{A\gamma(e\sigma(e))}:A\gamma(e\sigma(e))\to A\gamma(e\sigma(e)),\quad\theta_{l\m}=\gamma_l\m,\quad\theta_{m\m}=\sigma\gamma\m|_{A\gamma(e\sigma(e))}.
\end{align*}
Notice that $\gamma$ satisfies the condition (C1). Also $I_x$ is unital and $\gamma_{(x)}$ has globalization $\tilde{\gamma}_{(x)}$ given by $\tilde{\gamma}_{x}=\id_{A}$ and $\tilde{\gamma}_{g}=\sigma$. Thus, (C2) is also true. Moreover, $\gamma$ satisfies (C3) for $J_x=J_y=A$,
$\tilde{\gamma}_{\ta_x}=\tilde{\gamma}_x=\id_{A}$ and $\tilde{\gamma}_{\ta_y}=\gamma$.
By Proposition \ref{teo:globaliza1}, the action $\beta$ of $\G$ on $A$ given by
\[
{\beta}_{x}={\beta}_{y}=\id_A,\,\,\, {\beta}_{g}=\sigma,\,\,\, {\beta}_{h}=\gamma\sigma\gamma\m,\,\,\,{\beta}_{m}=\gamma\sigma,\,\,\, {\beta}_{m\m}=\sigma\gamma\m,\,\,\,{\beta}_{l}=\gamma,\,\,\,{\beta}_{l\m}=\gamma\m,
\]
is the globalization of the partial action $\theta$.}
\end{exe}

\section{Applications}

In this section we will show that, under appropriate assumptions, the Morita and Galois theories for a partial action of a connected groupoid $\G$ constructed via the functor $\Ext$ are strongly related with the Morita and Galois theories for the partial action of an isotropy group $\G(x)$. \vu

Throughout the rest of this paper we will assume that $\G$ is a connected finite groupoid, $x$ is a fixed object of $\G$, $\tau(x)$ is a transversal for $x$, $A$ is a ring and $\gamma=\big(I,\gamma_{\ta(x)},\gamma_{(x)}\big)$ is an object of $ \DA(A)$ that satisfies $I_{\ta\m_{y}}=I_x$ and $I_{\ta_{y}}=I_y$, for all $y\in \G_0$. We will also assume that $I_x$ is unital and $\gamma_{(x)}$ is a unital partial action, i.\,e., $I_h=I_x1_h$ where $1_h$ is a central idempotent of $I_x$ for all $h\in \G(x)$. Fix $\theta:=\Ext(\gamma)=(B_g,\theta_g\big)_{g\in\G}$  the partial action given by \eqref{afcon} and \eqref{dom}. \vu

\begin{lem1}\label{ob-partial-unital}  $B_g=A1_g$, where $1_g=\gamma_{\ta_{t(g)}}(1_{g_x})$ is a central idempotent of $A$ for all $g\in \G$.
\end{lem1}
\begin{proof}
Suppose that $I_x=A1_x$ and $I_h=I_x1_h$, where $1_x$ is a central idempotent of $A$ and $1_h$ is a central idempotent of $I_x$, for all $h\in \G(x)$. Then $a1_h=a(1_x1_h)=(a1_x)1_h=1_h(a1_x)=1_h(1_xa)=(1_h1_x)a=1_ha$, for all $a\in A$, $h\in \G(x)$.  Hence $1_h$ is a central idempotent of $A$ and $I_h=A1_h$, for all $h\in \G(x)$.
Using that $I_{\ta\m_{y}}=I_x$, $y\in \G_0$, we obtain from \eqref{dom} that $B_g=\gamma_{\ta_{t(g)}}(I_{g_x})$. Thus $B_g=A1_g$, where $1_g=\gamma_{\ta_{t(g)}}(1_{g_x})$ is a central idempotent of $A$, for all $g\in \G$.
\end{proof}

In what follows in this paper, we will assume that
$A=\oplus_{y\in\G_0}I_{y}=\oplus_{y\in\G_0}B_{y}$,
which is the appropriate context for working with invariant elements and trace maps.
\subsection{Invariant elements and trace map}

According to \cite{BP} an element $a\in A$ is called \emph{invariant} by $\theta$ if $\theta_g(a1_{g\m})=a1_g$, for all $g\in\G$. We will denote by $\A^\theta$ the set of all invariants of $A$. Clearly, $A^\theta$ is a subring of $A$. Similarly, $b\in  B_x^{\gamma_{(x)}}$ if and only if $\gamma_h(b1_{h\m})=b1_h$ for all $h\in \G(x)$. \vu

\begin{lem1}\label{inv}
	Let $b=\somay b_y\in A$.  Then $b\in A^\theta$ if and only if $b_x\in B_x^{\gamma_{(x)}}$ and $b_y=\gamma_{\ta_y}(b_x)$, for all $y\in\G_0$.
\end{lem1}

\begin{proof}
If $b\in A^\theta$ then $b_x1_h=b1_h=\theta_h(b1_{h\m})=\gamma_h(b_x1_{h\m})$ for all $h\in \G(x)$. Hence $b_x\in B_x^{\gamma_{(x)}}$.
Also, $\gamma_{\ta_y}(b_x)=\gamma_{\ta_y}\gamma_x\gamma\m_x(b_x1_x)=\theta_{\ta_y}(b1_{\ta\m_y})=b1_{\ta_y}=b_y$, for all $y\in \G_0$. Conversely, if $b_x\in B_x^{\gamma_{(x)}}$ and  $b=\somay\gamma_{\ta_y}(b_x)$ then \vu
\begin{align*}
	\theta_g(b1_{g\m})&=\gamma_{\ta_{t(g)}}\gamma_{g_x}\gamma\m_{\ta_{s(g)}}(b1_{g\m})\overset{\text{Lem.}\ref{ob-partial-unital}}{=}\gamma_{\ta_{t(g)}}\gamma_{g_x}\gamma\m_{\ta_{s(g)}}(b\gamma_{\ta_{s(g)}}(1_{g\m_x}))\\
	&=\gamma_{\ta_{t(g)}}\gamma_{g_x}\gamma\m_{\ta_{s(g)}}(\gamma_{\ta_{s(g)}}(b_x)\gamma_{\ta_{s(g)}}(1_{g\m_x}))=\gamma_{\ta_{t(g)}}\gamma_{g_x}(b_x1_{g\m_x})\\
	&=\gamma_{\ta_{t(g)}}(b_x1_{g_x})=b1_g,
	\end{align*}
for all $g\in\G$.\end{proof}

The \emph{trace map} $t_\theta:A\to A$ associated to the partial groupoid action $\theta$ is given by \[t_\theta(a)=\somag \theta_g(a1_{g\m}), \quad a\in A.\] It was proved in Lemma 4.2 of \cite{BP} that  $t_\theta(A)\subseteq A^\theta$ and consequently  $t_\theta:A\to A^\theta$. Similarly, the trace map $t_{\theta_{(x)}}:B_x\to B_x^{\theta_{(x)}}$ associated to the partial group action $\theta_{(x)}$ of $\G(x)$ on $B_x$ is given by $t_{\theta_{(x)}}(b)=\sum_{h\in\G(x)}\theta_h(b1_{h\m})$, $b\in B_x$. By Proposition \ref{sub} iii), $\theta_{(x)}=\gamma_{(x)}$ and whence $t_{\theta_{(x)}}=t_{\gamma_{(x)}}$.

\begin{lem1}\label{p52} Let $b_x\in B_x$ and $b_z=\gamma_{\ta_z}(b_x)\in B_z$. Then\vu \[t_{\theta}(b_z)=\sum_{y\in \G_0} \gamma_{\tau_y}(t_{\theta_{(x)}}(b_x)).\]
\end{lem1}
\begin{proof}
	Indeed,
	\begin{align*}
	t_{\theta}(b_z)&=\sum_{g\in \G}\theta_g(b_z1_{g\m})=\sum_{g\in \G}\gamma_{\tau_{t(g)}}\gamma_{g_x}\gamma\m_{\tau_{s(g)}}(\gamma_{\tau_{z}}(b_x)\gamma_{\tau_{s(g)}}(1_{g\m_x}))\\
	&=\sum_{s(g)=z}\gamma_{\tau_{t(g)}}\gamma_{g_x}(b_x1_{g\m_x})\overset{(l=g\ta_z)}{=}\sum_{s(l)=x}\gamma_{\tau_{t(l)}}\gamma_{l_x}(b_x1_{l\m_x}) \\
	&=\sum_{y\in \G_0}\sum_{l\in \G(x,y)}\gamma_{\tau_{y}}\gamma_{l_x}(b_x1_{l\m_x}) \overset{(\ta_yh=l)}{=}\sum_{y\in \G_0}\sum_{h\in \G(x)}\gamma_{\tau_{y}}\gamma_{(\ta_yh)_x}(b_x1_{(\ta_yh)\m_x})\\
	&=\sum_{y\in \G_0}\sum_{h\in\G(x)}\gamma_{\tau_{y}}\gamma_{h}(b_x1_{h\m})=\sum_{y\in \G_0}\gamma_{\tau_{y}}(\sum_{h\in \G(x)}\gamma_{h}(b_x1_{h\m}))\\
	&=\sum_{y\in \G_0}\gamma_{\tau_{y}}(t_{\theta_{(x)}}(b_x)).
	\end{align*}	
\end{proof}

\begin{lem1}\label{53}
Let $a=\sum_{z\in\G_0}b_z\in A$. Then \[t_\theta(a)=\sum_{y\in\G_0}\gamma_{\ta_y}(t_{\theta_{(x)}}(c_x)),\,\,\,\,\,\,\text{where } c_x:=\somaz \gamma\m_{\tau_z}(b_z)\in B_x.\]
\end{lem1}

\begin{proof}
	Indeed,
	\begin{align*}
	t_\theta(a)&=\somaz t_\theta(b_z)\overset{\text{Lem.}\ref{p52}}{=}\somaz\somay\gamma_{\ta_y}(t_{\theta_{(x)}}(\gamma_{\tau\m_z}(b_z)))\\
	&=\somay\gamma_{\ta_y}\left(t_{\theta_{(x)}}\left(\somaz \gamma_{\tau\m_z}(b_z) \right)\right)=\somay\gamma_{\ta_y}(t_{\theta_{(x)}}(c_x)).
	\end{align*}
\end{proof}

Now we will prove a result relating the surjectivity between $t_{\theta}$ and $t_{\theta_{(x)}}$ which will be useful later.

\begin{prop1} \label{t54}
	$t_\theta(A)=A^\theta$ if and only if $t_{\theta_{(x)}}(B_x)=B_x^{\theta_{(x)}}$.
\end{prop1}

\begin{proof}
	Assume that $t_\theta$ is onto. Take $b_x\in B_x^{\theta_{(x)}}$ and set $b=\somay\gamma_{\ta_y}(b_x)$ which lies in $A^\theta$ by Lemma \ref{inv}.
	Hence $b=t_\theta(a)$ for some $a\in A$. Since $A=\oplus_{z\in\G_0}B_z$, it follows
	that $a=\somaz b'_z$, with $b'_z\in B_z$ for each $z\in \G_0$. Then \[\somay\gamma_{\ta_y}(b_x)=b=t_\theta(a)\overset{\text{Lem.}\ref{53}}{=}\somay\gamma_{\ta_y}(t_{\theta_{(x)}}(c_x)),\]
	where $c_x=\somaz \gamma\m_{\tau_z}(b'_z)\in B_x$. Thus	$b_x=t_{\theta_{(x)}}(c_x)\in t_{\theta_{(x)}}(B_x)$. Conversely, notice that by Lemma \ref{inv},  each element in $A^\theta$ is of the form $b=\somay\gamma_{\ta_y}(b_x)$ with $b_x\in B_x^{\theta_{(x)}}$. By assumption, there exists $a_x\in B_x$ such that $t_{\theta_{(x)}} (a_x)=b_x$. Then, Lemma \ref{p52} implies that  $t_\theta(a_x)=\somay\gamma_{\ta_y}(t_{\theta_{(x)}}(a_x))=\somay\gamma_{\ta_y}(b_x)=b$.
\end{proof}

\subsection{Morita theory}

We start by recalling the definition of a Morita context. Given two unital rings $R$ and $S$,
bimodules $_RU_S$ and $_SV_R$, and bimodule maps $\mu: U\otimes_S V\to R$ and $\nu:V\otimes_R U\to S$,
the sixtuple $(R,S,U,V,\mu,\nu)$ is called a \emph{Morita context} (\emph{associated with $_RU$}) if the set
$$\left(
\begin{array}{cc}
R & U \\
V & S \\
\end{array}
\right)=\left\{\left(
\begin{array}{cc}
r & u \\
v & s \\
\end{array}
\right)\ \big|\ r\in R, s\in S, u\in U, v\in V\right\}
$$ with the usual addition and multiplication given by the rule
$$\left(
\begin{array}{cc}
r & u \\
v & s \\
\end{array}
\right)\left(
\begin{array}{cc}
r' & u' \\
v' & s' \\
\end{array}
\right)=\left(
\begin{array}{cc}
rr'+\mu(u\otimes v') & ru'+us' \\
vr'+sv' & \nu(v\otimes u')+ss' \\
\end{array}
\right)
$$ is a unital ring, which is equivalent to say that the maps $\mu$ and $\nu$ satisfy the
following two associativity conditions:
\begin{align}\label{associ-cond}
u\nu(v\otimes u')=\mu(u\otimes v)u'\quad \text{and}\quad v\mu(u\otimes v')=\nu(v\otimes u)v'.
\end{align}

We say that this context is \emph{strict} if the maps $\mu$ and $\nu$ are isomorphisms and, in this case,
the categories $_R\text{Mod}$ and $_S\text{Mod}$ are equivalent via the mutually inverse
equivalences $V\otimes_R-:\!_R\text{Mod}\to\!_S\text{Mod}$ and $U\otimes_S-:\!_S\text{Mod}\to\!_R\text{Mod}$. When it happens we also say that
the rings $R$ and $S$ are \emph{Morita equivalent}. The surjectivity  of  $\mu$ and $\nu$
is enough to ensure the strictness of the above context. Similar statements equally hold for right module categories. More details on Morita contexts can be seen, e.\,g., in $\S\,3.12$ of \cite{Jac}.

\subsubsection{Partial skew groupoid ring and partial skew group ring}
First we recall from \cite{BP} that the partial skew groupoid ring $R:=A\star_\theta\G$ associated to the partial groupoid action $\theta$ of $\G$ on $A$ is the set of all formal sums $\sum_{g\in \G}b_g\delta_g$, where $b_g\in B_g$, for each $g\in \G$. The addition in $R$ is the usual and the multiplication is given by
\[(b_g\delta_g)(b'_h\delta_h)=\left\{
\begin{array}{ll}
b_g\theta_g(b'_h1_{g\m})\delta_{gh}, & \text{if } \,\, (g,h)\in  \G_s\times\!_t\G, \\[.15cm]
0, & \text{otherwise.} \\
\end{array}\right.
\]
Similarly $S:=B_x\star_{\theta_{(x)}}\G(x)$ is the partial skew group ring associated to the partial group action $\theta_{(x)}$ of $\G(x)$ on $B_x$. Notice that $R$ and $S$ are unital associative rings, where $1_{R}=\sum_{y\in \G_0}1_y\delta_y$ and  $1_{S}=1_x\delta_x$.  In order to relate $R$ and $S$ we recall that
\begin{align}\label{for:ideal-action}
B_g=\gamma_{\ta_{t(g)}}(I_{g_x})\quad\text{and}\quad \theta_g(\gamma_{\ta_{s(g)}}(a))=\gamma_{\ta_{t(g)}}(\gamma_{g_x}(a)), \,\,\,\text{for all}\,\,\,g\in\G,\,\, a\in I_{g\m_x}.
\end{align}
In fact, the first identity of \eqref{for:ideal-action} can be seen in the proof of Lemma \ref{ob-partial-unital} and the second follows from \eqref{afcon}.

\begin{lem1}\label{lem-mc-skew} Let $R$ and $S$ be the rings defined above. Then
	\begin{enumerate}[\hspace{0.2cm}\rm i)]
		\item  $R1_S=\sum_{s(g)=x} B_g\delta_g$,\vspace*{.15cm}
		\item  $1_SR=\sum_{t(g)=x} B_g\delta_g$,\vspace*{.15cm}
		\item  $1_SR1_S=S$ and $R1_SR=R$.
	\end{enumerate}
\end{lem1}
\begin{proof} Let $g\in \G$ such that $s(g)=x$ and $a\in B_g$. By \eqref{for:ideal-action}, there exists $a'\in I_{g_x}$ such that $a=\gamma_{\ta_{t(g)}}(a')$ and whence
	\begin{align*}
	(a\delta_g)1_S&=(\gamma_{\ta_{t(g)}}(a')\delta_g)(1_x\delta_x)=\gamma_{\ta_{t(g)}}(a')\theta_g(\gamma_{\ta_{s(g)}}(1_{g\m_x}))\delta_{g}\\
	&\overset{\mathclap{\eqref{for:ideal-action}}}{=}\gamma_{\ta_{t(g)}}(a')\gamma_{\ta_{t(g)}}(\gamma_{g_x}(1_{g\m_x}))\delta_{g}=\gamma_{\ta_{t(g)}}(a'1_{g_x})\delta_{g}\\
	&=a\delta_g.
	\end{align*}
	For $g\in \G$ such that $s(g)\neq x$ and $a\in B_g$, we have $(a\delta_g)1_S=0$. Consequently, if $r=\sum_{g\in \G}a_g\delta_g\in R$ then $r1_S=\sum_{s(g)=x}a_g\delta_g$ and i) follows. The proof of ii) is similar.	

\vu

For iii), it is clear that $R1_SR\subset R$. For the reverse inclusion notice that
	\[\gamma_{\ta_{t(g)}}(a')\delta_g=\gamma_{\ta_{t(g)}}(a')\gamma_{\ta_{t(g)}}(1_{g_x})\delta_g=(\gamma_{\ta_{t(g)}}(a')\delta_{\ta_{t(g)}})(1_{g_x}\delta_{\ta\m_{t(g)}g}),\]
	for all $g\in \G$ and $a'\in I_{g_x}$. By ii), $1_{g_x}\delta_{\ta\m_{t(g)}g}\in 1_SR$. Hence $\gamma_{\ta_{t(g)}}(a')\delta_g\in R1_SR$. The equality $1_SR1_S=S$ follows in a similar way. 
\end{proof}

\vu

Fix the $(R,S)$-bimodule $U:=R1_S$ and the $(S,R)$-bimodule $V:=1_SR$. Define also the bimodule map $\mu: U\otimes_S V\to R$ (resp. $\nu:V\otimes_R U\to S$) given by $u\otimes v\mapsto uv$ (resp. $v\otimes u\mapsto vu$), for all $u\in U$, $v\in V$. With this notation we have the following result.

\begin{teo1}\label{teo:mc-skew} The sixtuple $(A\star_\theta\G, B_x\star_{\theta_{(x)}}\G(x),U, V,\mu,\nu)$ is a strict Morita context, i.\,e., $A\star_\theta\G$ and $B_x\star_{\theta_{(x)}}\G(x)$ are Morita equivalent.
\end{teo1}
\begin{proof}
	It is straightforward to check that $\mu,\nu$ satisfy \eqref{associ-cond}. By Lemma \ref{lem-mc-skew} iii), $\mu$ and $\nu$ are surjective maps.
\end{proof}

\subsubsection{Partial skew groupoid ring and subring of invariants}
Denote by $R:= A\star_{\theta}\G$ the partial skew groupoid ring. According to $\S\, 4$ of \cite{BP} $A$ is a  $(A^{\theta}, R)$-bimodule and a $(R, A^{\theta})$-bimodule. The left and the right actions of $A^\theta$ on $A$ are given by the multiplication of $A$, and the left (resp., right) action of $R$ on $A$ is given by $a\delta_g\cdot b=a\theta_{g}(b1_{g\m})$ (resp., $b\cdot a\delta_g=\theta_{g\m}(ba)$), for all $a\in A_g$, $b \in A$ and $g\in \G$. Moreover, the maps
\[\Gamma:A\otimes_{R} A\to A^{\theta}, \qquad a\otimes b\mapsto t_{\theta}(ab)\]
and
\[\Gamma':A\otimes_{A^{\theta}} A\to R, \qquad a\otimes b\mapsto \sum_{g\in \G}a\theta_g(b1_{g\m})\delta_g,\]
are bimodule morphisms. By Proposition 4.4 of \cite{BP}, the sixtuple $(R, A^{\theta}, A, A,\Gamma,\Gamma')$ is a Morita context. Furthermore, the sixtuple  $(B_x\star_{\theta_{(x)}} \G(x), B_x^{\theta_{(x)}}, B_x, B_x,\Gamma_x,\Gamma'_x)$ is also a Morita context constructed in a similar way as the previous one for the partial group action of $\G(x)$ on $B_x$ (see Theorem 1.5 of \cite{AL}).

\begin{cor1} \label{61}
	$\Gamma$ is surjective if and only if $\Gamma_x$ is surjective.
\end{cor1}
\begin{proof}
	Since $A$ (resp. $B_x$) is a unital ring, it is immediate to check that $\Gamma$ (resp. $\Gamma_x$) is surjective if and only if $t_{\theta}$ (resp. $t_{\theta_{(x)}}$) is surjective. Therefore the result follows from Proposition \ref{t54}.	
\end{proof}

\begin{lem1} \label{62}
	$\Gamma'$ is surjective if and only if $\Gamma'_{x}$ is surjective.
\end{lem1}
\begin{proof}
($\Rightarrow$) Let $b\delta_h\in B_x\star_{\theta_{(x)}}\G(x)$ with $b\in B_h\subseteq B_x$.  Consider $b\delta_h$ as an element of $R=A\star_\theta\G$. Since $\Gamma'$ is surjective, there exist $a_i,b_i\in A$, $1\leq i\leq r$, such that $\Gamma'\left(\sum_{1\leq i\leq r}a_i\otimes b_i\right)=b\delta_h$. By \eqref{for:ideal-action}, $B_{\ta\m_y}=\gamma_{\ta_x}(I_x)=\gamma_{x}(I_x)=I_x=B_x$ for all $y\in \G_0$. Similarly, $B_{\ta_y}=B_y$. Hence, $\theta$ is a group-type partial action. Thus, by Theorem 4.4 of \cite{BPP}, there are a global action $\beta$ of $\G_0^2$ on $A$ and a partial action $\mu$ of $\G(x)$ on $A\star_{\beta}\G_0^2$ such that
\[\varphi:A\star_\theta\G\to (A\star_{\beta}\G_0^2)\star_ \mu\G(x),\quad \varphi(a\delta_g)=a\delta_{(s(g),t(g))}\delta_{g_x},\quad a\in B_g, \,g\in \G,\]
is a ring isomorphism. Notice that $\varphi\circ\Gamma'\left(\sum_{1\leq i\leq r}a_i\otimes b_i\right)=\varphi(b\delta_h)=b\delta_{(x,x)}\delta_h$.
Since $A=\oplus_{y\in\G_0}B_{y}$, we can assume that $a_i=\somaz\gamma_{\ta_z}(a'_{i,z})$ and $b_i=\somaz\gamma_{\ta_z}(b'_{i.z})$ with $a'_{i,z},b'_{i,z}\in B_x$, for all $1\leq i\leq r$ and $z\in \G_0$. Then \vu
\begin{align*}
	\varphi\circ\Gamma'\left(\sum_{1\leq i\leq r}a_i\otimes b_i\right)&=\somai\somag a_i\theta_g(b_i1_{g\m})\delta_{(s(g),t(g))}\delta_{g_x}\\
	&=\somai\somag a_i\gamma_{\ta_{t(g)}}\gamma_{g_x}\gamma_{\ta\m_{s(g)}}(\gamma_{\ta_{s(g)}}(b'_{i,s(g)}1_{g\m_x}))\delta_{(s(g),t(g))}\delta_{g_x}\\
	&=\somai\somag a_i\gamma_{\ta_{t(g)}}\gamma_{g_x}(b'_{i,s(g)}1_{g\m_x})\delta_{(s(g),t(g))}\delta_{g_x}\\
	&=\somai\somag\gamma_{\ta_{t(g)}}(a'_{i,t(g)}\gamma_{g_x}(b'_{i,s(g)}1_{g\m_x}))\delta_{(s(g),t(g))}\delta_{g_x}.\
\end{align*} \vu
Hence, $$b\delta_{(x,x)}\delta_h=\somai\somag\gamma_{\ta_{t(g)}}(a'_{i,t(g)}\gamma_{g_x}(b'_{i,s(g)}1_{g\m_x}))\delta_{s(g),t(g))}\delta_{g_x}$$ which implies that
$$\somai a'_{i,x}\gamma_{h}(b'_{i,x}1_{h\m})=b\,\,\text{ and}\,\, \somai a'_{i,x}\gamma_l(b'_{i,x}1_{l\m})=0$$ if $l\neq h$, $l\in \G$.
Thus $\Gamma'_x(\somai a'_{i,x}\otimes b'_{i,x})=b\delta_h$ and $\Gamma'_x$ is surjective.
	
\vd
	
\noindent ($\Leftarrow$) Since $\varphi$ is a ring isomorphism, it is enough to check that for any element of the form $v= \gamma_{\ta_z}(a)\delta_{(y,z)}\delta_h$ in $(A\star_{\bt}\G_0^2)\star_ \mu\G(x)$,
	with $a\in A_h$ and $y,z\in\G_0$, there exists  $w\in A\otimes A$ such that $\varphi\circ\Gamma'(w)=v$. Observe that $a\delta_h\in B_x\star_{\theta_{(x)}}\G(x)$ and so there exist  $a_{x,i}, b_{x,i}\in B_x$ such that $a\delta_h =\Gamma'_x(\somai a_{x,i}\otimes b_{x,i})=\somal\somai a_{x,i}\gamma_l(b_{x,i}1_{l\m})\delta_l$. Thus
	\begin{align}\label{for-use}
	\somai a_{x,i}\gamma_h(b_{x,i}1_{h\m})=a \ \ \text{and} \ \ \somai a_{x,i}\gamma_l(b_{x,i}1_{l\m})=0 \ \ \text{if} \ l\neq h,\ l\in \G(x).
	\end{align}
	Now setting $a_i=\gamma_{\ta_z}(a_{x,i})$ and $b_i=\gamma_{\ta_y}(b_{x,i})$ in $A$ one has
	\begin{align*}
	\varphi\circ\Gamma'\left(\somai a_i\otimes b_i\right)&=\somai\somag a_i\theta_g(b_i1_{g\m})\delta_{(s(g),t(g))}\delta_{g_x}\\
	&=\somai\somag\gamma_{\ta_z}(a_{x,i})\gamma_{\ta_{t(g)}}\gamma_{g_x}\gamma\m_{\ta_{s(g)}}(\gamma_{\ta_y}(b_{x,i})\gamma_{\ta_{s(g)}}(1_{g\m_x})\delta_{(s(g),t(g))}\delta_{g_x}\\
	&=\somai\sum_{g\in\G(y,z)}\gamma_{\ta_z}(a_{x,i}\gamma_{g_x}(b_{x,i}1_{g\m_x}))\delta_{(y,z)}\delta_{g_x}\\
	&=\sum_{g\in\G(y,z)}\gamma_{\ta_z}\big(\somai a_{x,i}\gamma_{g_x}(b_{x,i}1_{g\m_x})\big)\delta_{(y,z)}\delta_{g_x}\\
	&\overset{\eqref{for-use}}{=}\gamma_{\ta_z}(a)\delta_{(y,z)}\delta_h.\
	\end{align*}
\end{proof}

Now we present the main result of this subsection which relates the strictness of the Morita contexts associate to $\theta$ (partial grupoid action) and to $\theta_{(x)}$ (partial group action).

\begin{teo1}\label{t63}
	If $A$ (resp., $B_x$) is a left $A^\theta$-module (resp., $B_x^{\theta_{(x)}}$-module) then the following statements are equivalent:
	\begin{enumerate}[\hspace{0.2cm}\rm i)]
		\item the Morita context  $(A^{\theta}, A\star_{\theta} \G, A, A,\Gamma,\Gamma')$ is strict,
		\item the Morita context  $(B_x^{\theta_{(x)}}, B_x\star_{\theta_{(x)}} \G(x),  B_x, B_x,\Gamma_x,\Gamma'_x)$ is strict.
	\end{enumerate}
\end{teo1}
\begin{proof}
	The result follows from Corollary \ref{61}	and Lemma \ref{62}.
\end{proof}
\subsection{Galois theory}
The notion of partial Galois extension for partial groupoid actions was introduced in \cite{BP} as a generalization of the classical one for group actions due to S.U. Chase, D.K. Harrison and A. Rosenberg in the global case \cite{CHR} and to M. Dokuchaev, M. Ferrero and the second named author in the partial case \cite{DFP}. This notion is closely related with the strictness of the Morita context between the partial skew groupoid ring and the subring of invariants.\vu

We say that $A$ is a \emph{$\theta$-partial Galois extension} of its subring of invariants $A^\theta$ if there exist elements $a_i,b_i\in A$, $1\leq i\leq r$, such that $\somai a_i\theta_g(b_i1_{g\m})=\delta_{y,g}1_y$, for all $y\in\G_0$, where $\delta_{y,g}$ denotes the Kronecker symbol. This is equivalent to say that the map
$\Gamma'$ defined in the previous subsection is a ring isomorphism (cf. Theorem 5.3 of \cite{BP}). Furthermore, in this case $A$, is also projective and finitely generated as a right $A^\theta$-module. (cf. assertion (ii) in  \cite[Theorem 5.3]{BP}).

\vu

The following theorem shows that the Galois theory for partial connected groupoid actions and that for partial group actions are indeed close.

\begin{teo1}\label{76}
	The following statements are equivalent:
\begin{enumerate}[\hspace{0.2cm}\rm i)]
		\item $A$ is a partial Galois extension of $A^\theta$ and $t_\theta$ is onto, \vu
		\item the Morita context $(A^{\theta}, A\star_{\theta} \G, A, A,\Gamma,\Gamma')$ is strict,\vu
		\item the Morita context $(B_x^{\theta_{(x)}}, B_x\star_{\theta_{(x)}} \G(x),  B_x, B_x,\Gamma_x,\Gamma'_x)$ is strict,\vu
		\item $B_x$ is a partial Galois extension of $B_x^{\theta_{(x)}}$ and $t_{\theta_{(x)}}$ is onto.
	\end{enumerate}
\end{teo1}
\begin{proof}
The result follows directly from  Proposition \ref{t54}, Theorem \ref{t63}, Theorem 5.3 of \cite{BP} and  Theorem 3.1 of \cite{AL}.
\end{proof}

\noindent \emph{Acknowledgements.} We thank the referee for his/her valuable comments and especially for suggesting the inclusion of Corollary 3.4.


\begin{thebibliography}{99999}
	
\bibitem[AL]{AL} J. \'Avila, J. Lazzarin, \emph{ A Morita context related to finite groups acting partially on a ring.}
Algebra Discrete Math. {\bf 3}, 49-60 (2009).

\bibitem[BFP]{BFP} D. Bagio, D. Fl\^ores and A. Paques, \emph{Partial actions of ordered groupoids on rings.} J. Algebra Appl. \textbf{9}, 501-517 (2010).

\bibitem[BP]{BP} D. Bagio and A. Paques, \emph{Partial groupoid actions: globalization, Morita theory, and Galois theory.} Comm. Algebra \textbf{40}, 3658-3678 (2012).

\bibitem[BPP]{BPP} D. Bagio, A. Paques and H. Pinedo, \emph{On partial skew groupoids rings.} arXiv:1905.12608.

\bibitem[BPi]{BPi} D. Bagio and H. Pinedo, \emph{On the separability of the partial skew groupoid ring.} \emph{S\~ao Paulo J. Math. Sci.} (2017) 11:370–384

\bibitem[CHR]{CHR} S.U. Chase, D. K. Harrison and A. Rosenberg, \emph{Galois theory and Galois cohomology of commutative rings.}  Mem. AMS \textbf{52}, 1-19 (1968).

\bibitem[DE]{DE} M. Dokuchaev and R. Exel, \emph{Associativity of crossed products by partial actions, enveloping actions and partial representations.} Trans. Amer. Math. Soc.
\textbf{357}, 1931-1952 (2005).


\bibitem[DFP]{DFP} M. Dokuchaev, M. Ferrero and A. Paques, \emph{Partial actions and Galois theory.}   J. Pure Appl. Algebra \textbf{208},
77-87 (2007).

\bibitem[E1]{EX} R. Exel, \emph{Circle actions on C*-algebras, partial automorphisms, and a generalized Pimsner-Voiculescu exact sequence.}
J. Funct. Anal. {\bf 122}, 361-401 (1994).

\bibitem[E2]{EX2} R. Exel, \emph{ Partial dynamical systems, fell bundles and applications.} Mathematical Surveys and Monographs, AMS  {\bf 224}, (2017).

\bibitem[Gi]{Gi} N. D. Gilbert, \emph{Actions and expansions of ordered groupoids.} J. Pure Appl. Algebra {\bf 198}, 175-195 (2005).

\bibitem[KL]{KL}J. Kellendonk and M. V. Lawson, \emph{Partial actions of groups.} Internat. J. Algebra Comput {\bf 14}, 87-114 (2004).

\bibitem[Ja]{Jac} N. Jacobson, \emph{Basic algebra II: Second Edition.}  Dover Publications, (2009).

\bibitem[ML]{ML} S. Mac Lane, \emph{Categories for the working mathematician}, Graduate Texts in
Mathematics, Springer-Verlag, Second edition, (1998).

\bibitem[NOP]{NYOP} P. Nystedt, J. \"{O}inert, and H. Pinedo, \emph{Artinian and noetherian partial skew groupoid rings.}  J. Algebra \textbf{503}, 433-452 (2018).

\bibitem[PT]{PT} A. Paques, and  T. Tamusiunas, \emph{The Galois correspondence theorem  for  groupoid actions.} J. Algebra, \textbf{509}, 105-123 (2018).
	
\end{thebibliography}
\end{document}